\documentclass[a4paper,twoside]{article}

\usepackage{a4}
\usepackage{amsmath,amssymb,amsthm}
\newcounter{tempref}

	{\begin{enumerate}
	} %
	{\end{enumerate}}
\newenvironment{alenumerate}%
	{\begin{enumerate}}%
	{\end{enumerate}}
	{\begin{enumerate}}%
	{\end{enumerate}}
  {\begin{pf*}{PROOF (#1).}}{\end{pf*}}
	{\setlength{\itemsep}{0pt}%
	\setlength{\topsep}{0pt}%
	\begin{trivlist}\item[]\begin{displaymath}%
	\begin{array}{rcll}}%
	{\end{array}\end{displaymath}\end{trivlist}}
	{\begin{list}{}{\setlength{\leftmargin}{0pt}%
			\setlength{\topsep}{1cm}}%
	\item{\bf Acknowledgements:}}%
	{\end{list}}

\newcommand{\rightqed} {\mystrut{0em}\hfill{$\Box$}}

\newcommand{\diam}{\operatorname{diam}}
\newcommand{\supdist}{\operatorname{supdist}}

\newcommand{\actson}{\cdot}

\newcommand{\Reals}{{\mathbb R}}

\newcommand{\argument}{\_\!\_}

\newtheorem{mythm}{Theorem}

\newtheorem{mylem}[mythm]{Lemma}

\newtheorem{mydefn}[mythm]{Definition}
\newtheorem{myrem}[mythm]{Remark}
\newtheorem{myexmp}[mythm]{Example}
\newtheorem{myobserve}[mythm]{Observation}

        {\begin{mydefn}\rm}{\end{mydefn}}
        {\begin{myexmp}\rm}{\end{myexmp}}
\newenvironment{rmrem}%
        {\begin{myrem}\rm}{\end{myrem}}
        {\begin{myprob}\rm}{\end{myprob}}
        {\begin{myobserve}\rm}{\end{myobserve}}

\renewcommand{\rightqed}{}

\newcommand{\eat}[1]{}

\begin{document}

\title{Linearizability of Non-expansive Semigroup Actions on Metric Spaces}
\author{Lutz Schr{\"o}der\\
DFKI-Lab Bremen and \\
Department of Mathematics and Computer Science,\\
 Universit\"at Bremen}
\date{}

\maketitle
\begin{abstract}
We show that a non-expansive action of a topological semigroup $S$ on a
metric space $X$ is linearizable iff its orbits are bounded. The crucial
point here is to prove that $X$ can be extended by adding a fixed point
of $S$, thus allowing application of a semigroup version of the
Arens-Eells linearization, iff the orbits of $S$ in $X$ are bounded.
\end{abstract}

\section*{Introduction}

By a well-known construction due to Arens and Eells~\cite{ArensEells56},
every metric space can be isometrically embedded as a closed metric
subspace of a normed (linear) space. Using this construction (or other
linear extensions like the free Banach space), one can
show~\cite{Pestov86,Megrelishvili06} that a non-expansive action $\pi$
of a topological semigroup $S$ on a metric space is linearizable, i.e.\
arises by restricting an action of $S$ by linear contractions on a
normed space $V$ to a metric subspace of $V$, if $\pi$ has a fixed point
$z$ (which then serves as the $0$ of $V$). The question of when an
action $\pi$ is linearizable in general thus reduces to the question of
when $\pi$ can be extended by adding a fixed point.

It is trivial to observe that if $X$ is bounded, then $\pi$ may be
extended by adding a fixed point: introduce a new point~$z$, make~$z$ a
fixed point of~$S$, and put $d(z,x)=c$ for all $x\in X$, where
$c>\diam(X)$. It is then easy to check that the distance function~$d$
thus defined on $X\cup\{z\}$ is a metric, and that the action of $S$ on
the extended space is non-expansive. Here, we improve on this
construction by giving a necessary and sufficient criterion: $\pi$ may
be extended by adding a fixed point iff its orbits are bounded sets. We
thus obtain an exact linearizability criterion: $\pi$ is linearizable
iff its orbits are bounded.

\section{Preliminaries}

Throughout the exposition, fix a topological semigroup $S$.  We shall
generally be concerned with \emph{non-expansive actions} $\pi:S\times
X\to X$, with $\pi(s,x)$ denoted as $s\actson x$, of $S$ on metric
spaces $(X,d)$, i.e.\ $d(s\actson x,s\actson y)\le d(x,y)$ for all $s\in
S$ and all $x,y\in X$. In the special case that $(X,d)$ is a real normed
space $V$, we say that $\pi$ is \emph{linear} if the translation maps
$\breve{s}:x\mapsto s\actson x$ are linear maps on $V$. In this case,
non-expansivity of $\pi$ means that the $\breve{s}$ are
\emph{contracting}, i.e.\ $\|s\actson x\|\le\|x\|$ for all $x$. We say
that a map $f:X\to Y$ is \emph{equivariant} w.r.t.\ an action of $S$ on
$Y$ if $f(s\actson x)=s\actson f(x)$ for all $x\in X$.

We note an observation from~\cite{Megrelishvili06}, omitting the
(straightforward) proof:
\begin{mylem}\label{lem:joint}
For a non-expansive action $\pi:S\times X\to X$ on a metric space
$(X,d)$, the following are equivalent.
\begin{enumerate}
\item The action $\pi$ is jointly continuous.
\item The action $\pi$ is separately continuous.
\item The restriction $\pi:S\times Y\to X$ to some dense subspace $Y$ of
  $X$ is separately continuous. 
\end{enumerate}
\end{mylem}
We shall henceforth implicitly include the requirement that $S\times
X\to X$ is continuous in the term \emph{non-expansive action} (thus
avoiding the term `non-expansive continuous action', which is a bit of a
mouthful). As an immediate consequence of the preceding lemma, we obtain
the following extension result~\cite{Megrelishvili06}:
\begin{mylem}\label{lem:banach}
A linear non-expansive  action of $S$ on a normed space $V$
extends (uniquely) to a linear non-expansive  action of $S$ on
the completion of $V$.
\end{mylem}
We denote the orbit $\{s\actson x\mid s\in S\}$ of $x\in X$ under $S$ by
$S\actson x$. Note that orbits need not be disjoint, elements of an
orbit need not have the same orbit, and $x$ need not be contained in its
orbit $S\actson x$. In case $S$ is a monoid, however, $x\in S\actson x$
for all $x\in X$.

\section{Fixed Points and Linearizations}

We now give the announced criterion for extendibility by a fixed point:

\begin{mythm}\label{thm:main}
Let~$(X,d)$ be a metric space equipped with a non-expansive 
action of $S$. Then the following are equivalent:
\begin{enumerate}
\item The space $X$ can be extended by adding a fixed point of $S$,
  i.e.\ there exists a metric space~$(Y,d)$ equipped with a
  non-expansive action of $S$ that has a fixed point, and an isometric
  and equivariant embedding of~$X$ into~$Y$.
\item The orbits $S\actson x$ of~$S$ in~$X$ are bounded sets.
\end{enumerate}
\end{mythm}
The following definition will be useful in the proof:
\begin{mydefn}
Let $(X,d)$ be a metric space. For $A\subseteq X$ and $x\in X$, we put
\begin{equation*}
\supdist(x,A)=\sup_{y\in A}\, d(x,y)\in[0,\infty].
\end{equation*}
\end{mydefn}

\begin{proof}(Theorem~\ref{thm:main})
\emph{(1)$\implies$(2):} We can assume  that $X$ is a subspace
of~$Y$. Let $z\in Y$ be a fixed point of~$S$. Then we have, for $x\in X$
and $s,t\in S$,
\begin{align*}
d(s\actson x,t\actson x)  & \le d(s\actson x,z) + d(z,t\actson x) \\
& = d(s\actson x,s\actson z) + d(t\actson z,t\actson x) \\
& \le 2d(x,z),
\end{align*}
i.e.\ $\diam(S\actson x)\le 2d(x,z)$.

\emph{(2)$\implies$(1):} To begin, we reduce to the case that $S$ is a
monoid, as follows. For a semigroup $S$, we have the free monoid $S_e$
over $S$, constructed by taking $S_e=S\cup\{e\}$, where $e\notin S$, and
putting $es=se=s$ for all $s\in S_e$. The action of $S$ on $X$ is
extended to a (non-expansive) action of $S_e$ by $e\actson x=x$ for all
$x\in X$. The orbits $S_e\actson x=\{x\}\cup S\actson x$ are bounded (by
$d(x,s\actson x)+\diam(S\actson x)$ for any $s\in S$). By the monoid
case of the theorem, we obtain an extended space $(Y,d)$ in which $S_e$
has a fixed point $z$; the action of $S_e$ restricts to an action of $S$
on $Y$, and $z$ trivially remains a fixed point of $S$.

When $S$ is a monoid, then $x\in S\actson x$ for all $x\in X$.  We can
assume w.l.o.g.\ that there exists a point $x_0\in X$ which is not fixed
under $S$. We put $Y=X\cup\{z\}$, where $z\notin X$,
\begin{equation*}\textstyle
d(z,x)=d(x,z) = \supdist (x_0,S\actson x)
\end{equation*}
for $x\in X$, and $d(z,z)=0$. We have to check that this makes $(Y,d)$ a
metric space. To begin, $d(x,z)>0$ for $x\in X$: we have $\supdist
(x_0,S\actson x_0)>0$ because $x_0$ is not fixed under $S$, and for
$x\neq x_0$, $\supdist(x_0,S\actson x)\ge d(x_0,x)>0$ (using $x\in
S\actson x$).  Symmetry holds by construction. Moreover, for $x\in X$,
$d(x_0,s\actson x)\le d(x_0,x)+d(x,s\actson x) \le d(x_0,x) +
\diam(S\actson x)$ for all $s\in S$ (again using $x\in S\actson x$) and
hence $d(x,z)\le d(x_0,x) + \diam(S\actson x) <\infty$ by~(2).  It
remains to prove the triangle inequality. There are only two non-trivial
cases to prove:
\begin{alenumerate}
\item $d(x,z) \le d(x,y) + d(y,z)$ for $x,y\in X$, and
\item $d(x,y) \le d(x,z) + d(y,z)$ for $x,y\in X$.
\end{alenumerate}

\emph{Ad (a):} Let $s\in S$. Then $d(x_0,s\actson x)\le d(x_0,s\actson
y)+d(s\actson y,s\actson x)\le d(x_0,s\actson y)+d(y, x)$. Thus,
$\supdist(x_0,S\actson x)\le d(x,y) + \supdist(x_0,S\actson y)$.

\emph{Ad (b):} We have
\begin{align*}
d(x,y) 	&\le d(x,x_0) + d(y,x_0)\\
	&\le \supdist(S\actson x,x_0) + \supdist(S\actson y,x_0)\\
	&= d(x,z) + d(y,z),
\end{align*}
where the second inequality uses $x\in S\actson x$.

We then extend the action of $S$ to $Y$ by letting $z$ be fixed under
$S$. It is clear that this really defines an action of $S$; we have
to check that this action is non-expansive. For $x\in X$ and $s\in S$,
we have
\begin{align*}
d(s\actson x,s\actson z) & = d(s\actson x,z) \\
 & = \supdist(x_0,S\actson (s\actson x))\\
 & \le \supdist(x_0,S\actson x)\\
 & = d(x,z),
\end{align*}
where the inequality uses $S\actson (s\actson x)  \subseteq S\actson x$.

It remains to prove that $S\times Y\to Y$ is continuous, i.e.\ by
Lemma~\ref{lem:joint} that the orbit maps $S\to Y,s\mapsto s\actson y$,
are continuous. For $y\in X$, this follows from continuity of the action
on $X$, and for $y=z$, the orbit map is constant.  
\rightqed
\end{proof}

\begin{rmrem}
In case $S$ is a group, one can identify the space $Y$ constructed in
the above proof with the subspace $\{\{x\}\mid x\in X\}\cup\{S\actson
x_0\}$ of the space of bounded subsets of $X$, equipped with the
Hausdorff pesudometric
\begin{equation*}
d(A,B):=\max\{\sup\{d(a,B):   a \in A\},  \sup\{d(A,b):   b \in B\}\}
\end{equation*}
and the natural action taking $A$ to $s\actson A$ for $s\in S$.  For
arbitrary semigroups $S$, however, orbits will in general fail to be
fixed points under the natural action.
\end{rmrem}

\begin{rmrem}
Analogously as in the proof of (1)$\implies$(2) in the above theorem,
one shows that for $x,y\in X$, $\diam(S\actson y)\le
2d(x,y)+\diam(S\actson x)$. Thus, for boundedness of all orbits it
suffices to require that there exists a bounded orbit.
\end{rmrem}

We now briefly recall the Arens-Eells extension of a pointed metric
space $(X,d,z)$ (i.e. $z\in X$). One constructs a real normed space
$(A(X),\|\argument\|)$ by taking as the elements of $A(X)$ the formal
linear combinations
\begin{equation*}
\sum_{i=1}^n c_i(x_i-y_i),
\end{equation*}
with $x_i,y_i\in X$ and $c_i\in\Reals$ and putting for $u\in A(X)$
\begin{equation*}
\|u\|=\inf\left\{\sum_{i=1}^n |c_i|d(x_i,y_i)\mid u=\sum_{i=1}^n
c_i(x_i-y_i)\right\}.
\end{equation*}
The space $(X,d)$ is isometrically embedded into $A(X)$ (as a closed
subspace) by taking $x\in X$ to $x-z$. It is shown
in~\cite{Megrelishvili06} (Proposition 2.10) that a non-expansive 
action of $S$ on $X$ can be extended to a linear non-expansive
 action of $S$ on $A(X)$ by putting
\begin{equation*}
s\actson \sum_{i=1}^n c_i(x_i-y_i) = \sum_{i=1}^n c_i(s\actson
x_i-s\actson y_i).
\end{equation*}
(A similar construction can be found already in~\cite{Pestov86};
moreover, the Arens-Eells extension may be replaced by other linear
extensions~\cite{Flood84}, e.g.\ the free Banach space over $X$ as
in~\cite{Pestov86}.)

We then immediately obtain the announced exact linearizability
criterion. 
\begin{mythm}\label{thm:linear}
For a non-expansive  action of $S$ on a metric space~$(X,d)$,
the following are equivalent:
\begin{enumerate}
\item There exists a Banach space~$V$, equipped with a linear
  non-expansive  action of~$S$, and an equivariant isometric
  embedding of~$(X,d)$ into~$V$.
\item The orbits~$S\actson x$ of~$S$ in~$X$ are bounded sets.
\end{enumerate}
\end{mythm}
\begin{proof}
\emph{(1)$\implies$(2):} By the corresponding direction of
Theorem~\ref{thm:main}, as $0$ is a fixed point of $S$ in $V$.

\emph{(2)$\implies$(1):} By Theorem~\ref{thm:main}, we may assume that
$S$ has a fixed point $z$ in $X$. By Lemma~\ref{lem:banach}, it suffices
to construct $V$ as a normed space. We thus may take $V$ as the
Arens-Eells extension of $(X,d,z)$, equipped with the $S$-action
described above.  \rightqed
\end{proof}

 \bibliographystyle{abbrv}
\bibliography{fixedpoint}

\end{document}